\documentclass[11pt]{article}
\usepackage{amsfonts,color}
\usepackage{amssymb}
\usepackage{amsmath}
\usepackage{amsthm}
\usepackage[round]{natbib}
\usepackage{enumerate}
\usepackage{graphicx}
\usepackage{wrapfig}
\usepackage{rotating}
\usepackage{epsfig}
\usepackage[left=1in, right=1in, top=1in, bottom=1in]{geometry}
\usepackage{mdwlist}
\usepackage{secdot}
\usepackage{comment}
\usepackage{mathrsfs}
\usepackage{mdwlist}
\usepackage{multirow}
\usepackage{lscape}
\usepackage{pdflscape}
\usepackage{array}
\usepackage[compact]{titlesec} 
\usepackage{times}
\usepackage{tikz}
\usepackage{subfigure}
\usepackage[justification=centering]{caption}
\usepackage{setspace}
\usepackage{color}
\usepackage{authblk}
\usepackage{wrapfig}
\usepackage{hyperref}
\usepackage{makeidx}
\usepackage[title]{appendix}
\usepackage[ruled]{algorithm2e}
\usepackage{graphicx}
\usepackage{float}

\graphicspath{ {./images/} }

\bibpunct{[}{]}{,}{n}{}{;}
\titlespacing*{\section}{0pt}{*2}{*0}
\titlespacing*{\subsection}{0pt}{*2}{*0}
\titlespacing*{\subsubsection}{0pt}{*1}{*0}

\newtheoremstyle{compacttheorem}
{3pt}
{3pt}
{}
{}
{}
{}
{}
{}

\renewcommand{\arraystretch}{1.12}

\title{Incorporating Equity into the School Bus Scheduling Problem}
\author[a]{Dipayan Banerjee}
\author[b]{Karen Smilowitz}
\affil[a, b]{ Department of Industrial Engineering and Management Sciences, Northwestern University, 2145 Sheridan Road, Evanston, IL 60208, United States \vspace{2.5ex}}
\affil[a]{ dipayanbanerjee2019@u.northwestern.edu}
\affil[b]{ ksmilowitz@northwestern.edu}

\begin{document}
\maketitle

\vspace{-0.9cm}

\subsubsection*{\centering Abstract}

\indent We consider the school bus scheduling problem (SBSP) which simultaneously determines school bell times and route schedules. Often, the goal of the SBSP is to minimize the number of buses required by a school district. We extend a time-indexed integer programming model to incorporate additional considerations related to equity and efficiency. We seek to equitably reduce the disutilities associated with changing school start times via a minimax model, then propose a lexicographic minimax approach to improve minimax solutions. We apply our models to randomized instances based on a moderately-sized public school district to show the impact of incorporating equity.

\vspace{0.1cm}

\noindent \textbf{Keywords:} school bus scheduling, integer programming, equity, fairness, lexicographic minimax 

\vspace{-0.2cm}

\section{Introduction}

Since the economic recession in 2008, public school districts across the United States have faced budget crises as stagnant funding levels have often been insufficient to offset rising operational costs. Between 2008 and 2015, state funding per student for public schools declined in 29 states; in some states, state funding in 2017-18 had yet to return to pre-2008 levels (Leachman et al., \citeyear{punishing}).  Often, budget cuts entail a reduction in transportation spending (see recent examples in Kentucky (Spears, \citeyear{kentuckycuts}) and Missouri (Rowe, \citeyear{missouricuts})), necessitating more cost-efficient use of available transportation resources.  

Public school transportation spending in the United States for 2012-13, the last school year with full data, totaled about \$23 billion - approximately 4 percent of total public education spending (Urban Institute, \citeyear{national_costs}). Transportation costs for public school systems often include a component proportional to the total number of buses used. A typical contract for outsourced transportation services may include a fixed per-bus cost structure, with additional costs - significantly lower than the fixed cost - for multiple uses of the same bus. Previous analyses have encountered fixed costs per bus in Drummondville, Qu{\'e}bec (Desrosiers et al., \citeyear{transcol}) and New Haven, Connecticut (Swersey and Ballard, \citeyear{swersey}), among other places. One way to reduce the number of buses required by a school district, thereby reducing the associated outsourcing costs (or capital costs if the district owns the buses), is to schedule schools and routes in such a way that buses can be reused to serve multiple routes. Consider an example with two schools, each with two morning bus routes. In Figure 1(a), both schools start at the same time. Because all four routes run simultaneously, four buses are required to complete the routes. In Figure 1(b), School A starts 45 minutes before School B. Only two buses are required to complete the routes: one bus can service routes A1 and B1, while the other bus can service A2 and B2.

\begin{figure}[h]
	\includegraphics[width=\textwidth]{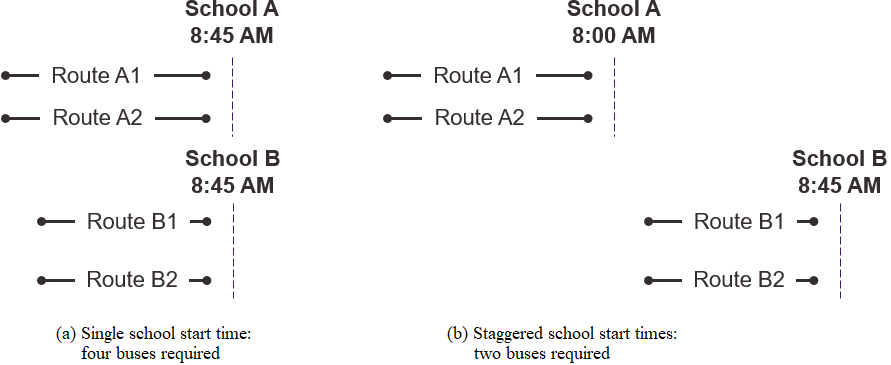} 
	\centering
	\caption{School and route scheduling to minimize the number of buses required}
\end{figure}
\vspace{-0.05in}

The combined problem of route scheduling and school bell time scheduling, with the primary objective of minimizing the number of buses required, has been studied in operations research as the \textit{school bus scheduling problem}, or SBSP (F{\"u}genschuh, \citeyear{fugenschuh}; Raff, \citeyear{raff}). The SBSP is a part of the more general \textit{school bus routing problem} (SBRP), which consists of five components: data preparation, bus stop selection, bus route generation, school bell time adjustment, and route scheduling (Park and Kim, \citeyear{park_kim_review}).

Our collaboration with a moderately-sized public school district has motivated a need to explicitly consider equity in the SBSP. Changing school start times can have both positive and negative effects on communities. For example, studies suggest that later school start times are correlated with improved behavior and academic performance in adolescents (Owens et al., \citeyear{owens}; Wolfson et al., \citeyear{wolfson}). Based on this research, an unsuccessful bill was recently proposed in California which would have required middle and high schools to start at 8:30 a.m.\ or later (Racker, \citeyear{california}). At the same time, some parents express concerns that later school start times interfere negatively with after-school activities (Dunietz et al., \citeyear{dunietz_start_times}). An optimal schedule should ensure that the effects of any school time changes are distributed equitably across the district. Specifically, our partner school district prefers schedules with lower absolute deviation in start times from the current schedule, and the district prefers any such absolute deviations to be distributed equitably.

The remainder of the paper is organized as follows. Section 2 reviews the relevant literature. Section 3 details the base integer programming model formulation for the SBSP and our approach to equity. Section 4 presents methods for integrating equity in regard to school start time changes. Section 5 contains numerical studies of our models. Finally, Section 6 contains our concluding remarks.
 
\section{Literature Review}
A thorough survey on mathematical approaches to the SBRP (through 2009) can be found in a review by Park and Kim (\citeyear{park_kim_review}). While numerous recent studies (e.g., B{\"o}gl et al., \citeyear{bogl}; Chen et al., \citeyear{chen}; Hashi et al., \citeyear{hashi}; Kim et al., \citeyear{kim_kim_park}; Shafahi et al., \citeyear{shafahi}; Yan et al., \citeyear{yan}) discuss ``school bus scheduling" in various contexts, these papers focus on route scheduling without considering school bell time adjustment. There are limited studies which examine school bell time adjustment in conjunction with route scheduling (i.e., the SBSP).

One such study is Desrosiers et al.\ (\citeyear{transcol}) which proposes a framework to solve parts of the SBRP sequentially: (1) route generation, followed by (2) school bell time scheduling, followed by (3) route scheduling. As such, the route generation component (1) provides the required route data for the SBSP components (2) and (3). Desrosiers et al.\ (\citeyear{transcol}) establish and use the following proposition, restated formally by Zeng et al.\ (\citeyear{zeng_sbsp}):

\begin{quote}
	\textbf{Proposition 1:} \textit{Given a set of routes with fixed starting and ending times, the minimum number of buses required to complete all the routes is equal to the maximum number of routes in operation in the same time period.}
\end{quote}

\noindent In practice, a bus which completes a route requires a transition time to travel to the starting location of its next route. As a result, the actual number of buses required will likely be more than the number implied by Proposition 1. However, if the transition time is included in the operational time of each route (i.e., a route is ``in operation" during the transition period), then Proposition 1 can be used to accurately determine the minimum number of buses required. Many approaches to route scheduling (e.g., Bookbinder and Edwards, \citeyear{bookbinder_SBSP}; F{\"u}genschuh, \citeyear{fugenschuh}; Spada et al., \citeyear{spada}) use exact travel times between routes, explicitly considering the travel distance between the ending point of one route and the starting point of the next route serviced by the same bus. As discussed in Section 3, we assume that the transition time between routes is constant (e.g., 10 or 15 minutes) and is added to the true duration of the route, as in Zeng et al.\ (\citeyear{zeng_sbsp}). Desrosiers et al.\ (\citeyear{transcol}) formulate the problem of scheduling school bell times as an integer linear program (ILP), which they solve using decomposition and heuristic approaches. Route schedules are determined heuristically via iteratively solving transportation problems. Their approach for the entire SBRP found significant savings when applied to the Saint-Fran\c{c}ois school board in Drummondville, Qu{\'e}bec, comprised of 47 schools and over 17,000 students.

F{\"u}genschuh and Martin (\citeyear{fugenschuh2006}) present an ILP formulation to model both components of the SBSP concurrently, which they solve using a two-stage heuristic to minimize seven lexicographically ordered objectives: number of buses required, total driving time without passengers, waiting time of buses between trips, total waiting time of students at school, total waiting time of students at transfer bus stops, total absolute change of school start times, and total absolute change of route start times. F{\"u}genschuh (\citeyear{fugenschuh}) solves a similar model to minimize the first two objectives via branch-and-cut methods. When applied to five school districts in Germany, this model found that the number of buses required could be reduced by 10-25\%. F{\"u}genschuh (\citeyear{fugenschuh2011}) later improves the model's lower bounds using a set partitioning reformulation. 

Using Proposition 1, Zeng et al.\ (\citeyear{zeng_sbsp}) present a time-indexed ILP model for the SBSP based on machine scheduling and bin packing problems. While Desrosiers et al.\ (\citeyear{transcol}) and F{\"u}genschuh (\citeyear{fugenschuh}) determine a sequence of routes for each bus, Zeng et al.\ (\citeyear{zeng_sbsp}) focus solely on determining route arrival times and school start times, noting that the associated problem of assigning buses to routes given a route schedule can be solved separately in polynomial time using an algorithm detailed by Olariu (\citeyear{olariu}). For cases where all routes associated with a school must arrive at the same time, a LP relaxation-based rounding algorithm with a probabilistic Chernoff upper bound and a greedy algorithm are presented. An LP relaxation-based rounding algorithm is presented for a slightly modified formulation of the general SBSP without the condition on time window lengths, and numerical tests are conducted on randomized instances of varying sizes. Park and Kim (\citeyear{park_kim_review}) discuss the need for a generalized SBRP; in this paper we demonstrate the potential of the model proposed by Zeng et al.\ (\citeyear{zeng_sbsp}) to serve as a generalized SBSP upon which modifications can be made to meet the planning needs of a particular school system.

Equity, in the context of transportation systems, ``refers to the distribution of impacts (benefits and costs) and whether that distribution is considered fair and appropriate" (Litman, \citeyear{litman_equity}). Savas (\citeyear{savas_equity}) discusses four broad characterizations of equity in the provision of public services: (1) equal payment, (2) equal output results or metrics, (3) equal provided inputs, and (4) equal satisfaction of demand. Of these, our approach focuses on equal output metrics, e.g., how changes in start times impact the student population. Equal payment is largely irrelevant to public education as all students have the same educational rights regardless of their ability to pay. Equal inputs and equal satisfaction of demand are accomplished outside of our analysis, as school districts should provide an adequate number of schools and all eligible students are provided transportation. Marsh and Schilling (\citeyear{marsh_equity}) introduce a framework for incorporating equity in operations research models (with a focus on facility location) which we employ in our paper. First, they define the \textit{determinants} of equity - the effects which should be distributed fairly. Next, they define the \textit{groups} across which equity is to be pursued. Finally, they review 20 potential \textit{measures} of inequity to be minimized.

In the SBRP literature, a limited number of studies pursue equity as an objective. In selecting bus stops and constructing routes for a single school, Bowerman et al.\ (\citeyear{bowerman}) and Li and Fu (\citeyear{li_SBRPequity}) aim to balance the number of students on each route and the length of each route. To achieve this, Bowerman et al.\ (\citeyear{bowerman}) include nonlinear measures of inequity in their objective function, while Li and Fu (\citeyear{li_SBRPequity}) incorporate inequity minimization into a heuristic improvement strategy. Delgado and Pacheco (\citeyear{delgado_SBRPequity}) minimize the cost of the costliest route in a modification to the classical vehicle routing problem (VRP), and Pacheco et al.\ (\citeyear{pacheco_SBRPequity}) simultaneously minimize both the maximum route length and the total route length in another modification to the VRP; both models are applied to the construction of school bus routes. Spada et al.\ (\citeyear{spada}) aim to balance students' \textit{time loss} by minimizing the maximum time loss across all students as a secondary objective. Time loss is the sum of a student's \textit{delay} (the difference in time between the ``actual journey time [between the student's home and school] and the shortest possible time between home and school") and \textit{waiting time} (the difference between the student's arrival time at school and the school's start time). Equity is one of the objectives considered in a recent study motivated by potential start time changes to Boston Public Schools (Bertsimas et al., \citeyear{mit_sbsp}): minimizing the variance (across neighborhoods or communities) in satisfaction with a potential new schedule is included as an objective in a quadratic optimization model.

Equity-based objectives also appear in recent operations research approaches to other problems in public education and public service at large.  Bouzarth et al.\ (\citeyear{bouzarth_districting}) minimize socioeconomic variation across schools in conjunction with transportation costs when assigning students to public schools, while Campbell et al.\ (\citeyear{campbell_equity}) minimize the maximum arrival time of relief supplies after a disaster. 

Further examples of equity considerations in the community-based operations research literature can be found in recent reviews by Leclerc et al.\ (\citeyear{leclerc}) and Johnson et al.\ (\citeyear{johnson}). In this paper, we integrate equity considerations into the SBSP as has been done in other public sector operations research problems. Our primary contributions to the literature are threefold. Firstly, we develop minimax and lexicographic minimax ILP models which simultaneously promote equity and efficiency in school bell time and school bus schedules with respect to disutilities arising from changing school start times. Through our partnership with an urban public school district, we ensure that our models, which are guided by the familiar idea of improving the condition of those who are worst-off, are readily accessible to policymakers. Secondly, we develop an iterative algorithm for solving the lexicographic minimax model which avoids numerical issues which may lead to incorrect solutions. Finally, we demonstrate that the lexicographic model produces significantly improved schedules while maintaining computational tractability. As such, the lexicographic approach presented in this paper serves as a practical framework for public school districts which seek to design equitable school and bus schedules, and as a potential foundation for operations researchers who seek to incorporate equity considerations into other modern transportation and logistics problems.

\section{Preliminaries}

In this section, we present the base model of school bus scheduling and a general approach to incorporate equity.

\subsection{Modeling school bus scheduling}
The formal definition of the problem, adapted from Zeng et al.\ (\citeyear{zeng_sbsp}), is as follows:
\vspace{0.1in}

\noindent \textbf{School Bus Scheduling Problem (SBSP):} \textit{Given a set of schools, each with a set of morning and afternoon routes, time windows for school start and end times, and associated time windows for route arrival and departure, determine the starting and ending time of each school and the arrival and departure times for all morning and afternoon routes, respectively, such that time windows are satisfied, and assign routes to buses. The objective is to minimize the number of buses required to complete all of the routes.} 
\vspace{0.1in}

The problem statement is necessarily broad, as each school district faces a unique problem setting in practice. In the literature, additional assumptions are included based on the motivating case. Our modeling assumptions, which reflect the problem setting of our partner school district, are as follows:

\begin{itemize}
	\item[(i)] School day lengths are sufficiently large such that all morning routes are completed before afternoon routes begin.
	\item[(ii)] The district area is geographically compact such that the travel times between the ending point of one route and the starting point of another route are modeled as uniform and short enough so that reusing buses on multiple routes is practical.
	\item[(iii)] Each route is associated with a single school (i.e., no mixed loads).
\end{itemize}

\noindent The second assumption implies that if the constant travel time between routes is added to the route length, Proposition 1 can be used to determine the minimum number of buses required instead of directly assigning routes to buses. The constant transition time may also provide an additional buffer in case of delays, increasing the reliability of the transportation system. For conciseness, we present the part of our ILP model associated with morning scheduling. In Appendix A, we show how the extension to include afternoon scheduling can be implemented easily. 

\subsubsection{ILP Formulation}
The morning time interval during which all routes must arrive and all schools must start is discretized over a total of $M$ periods. For example, the interval 8:00 AM to 8:15 AM discretized by 5 minutes gives $M = 4$, where the times are \{8:00, 8:05, 8:10, 8:15\}. We are given $N$ schools, where the $n^{th}$ school is associated with $\Gamma_n$ morning routes. The $i^{th}$ route associated with school $n$ has discretized length (duration) of $r_{i, n}$ consecutive time periods, including a constant transition time added to each route. The parameter $\alpha$ denotes the maximum number of time periods before the school start time that a route can arrive, and $\beta$ denotes the minimum number of time periods before the start time that a route must arrive. Figure 2 illustrates a scenario where time is discretized by five minute intervals; if $\alpha = 4$, $\beta = 2$, and a school starts at 9:00 AM, then the routes associated with the school may arrive at 8:40 AM, 8:45 AM, or 8:50 AM.

\begin{figure}[h]
	\includegraphics{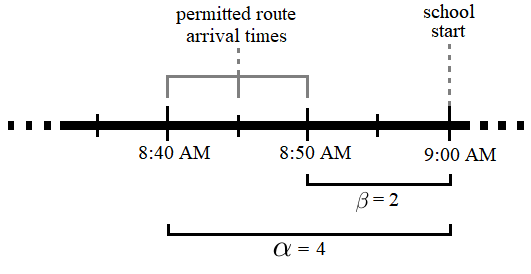} 
	\centering
	\caption{Arrival time window example}
\end{figure}

Define binary variables $x_{i, n}^{(m)}$ for all $n \in [N], i \in \Gamma_{n}, m \in [M]$ and $u_{n}^{(m)}$ for all $n \in [N], m \in \{\beta+1, \ldots, M\}$ as such: $x_{i, n}^{(m)} = 1$ if route $i \in \Gamma_{n}$ arrives at time $m$ and 0 otherwise; $u_{n}^{(m)} = 1$ if school $n$ starts at $m$ and 0 otherwise. Note that including school start time variables is not necessary in the base formulation to find the optimal number of buses; it is easy to show that start times can be inferred if the problem is formulated using only $x_{i, n}^{(m)}$ and $z$ variables (as in Zeng et al., \citeyear{zeng_sbsp}) given the time window parameters. However, start time variables are required for our model extensions. The use of binary variables to track school start times is similar to the approach used by Desrosiers et al.\ (\citeyear{transcol}) but is a departure from F{\"u}genschuh's (\citeyear{fugenschuh}) more recent approach of using an integer variable to denote each school's start time. 

Define $z$ as the minimum number of buses required to service all routes. Our goal is to determine a route and school bell schedule, defined by the values of the binary variables, such that $z$ is minimized. We formulate the SBSP as follows:
\vspace{-0.9cm}

\begin{align}
&min~~ z \nonumber \tag{ILP 1}\\
s.t.~~& \sum_{m=\beta+1}^{M}u_{n}^{(m)} = 1 &\forall n \in [N] \tag{1.1} \\
& u_{n}^{(m)} \leq \sum_{t=\max \{m-\alpha,1\}}^{m-\beta}x_{i,n}^{(t)} &\forall n \in [N], i \in [\Gamma_{n}], m \in \{\beta+1, \ldots, M\} \tag{1.2}\\ 
&\sum_{n=1}^{N}\sum_{i=1}^{\Gamma_{n}}\sum_{t=m}^{\min\{m+r_{i, n}-1, M \}}x_{i,n}^{(t)}\leq z &\forall m \in [M] \tag{1.3}\\
&x_{i,n}^{(m)} \in \{0,1\} &\forall n \in [N], i \in [\Gamma_{n}], m \in [M] \tag{1.4} \\
&u_{n}^{(m)} \in \{0,1\} &\forall n \in [N], m \in \{\beta+1, \ldots, M\} \tag{1.5}
\end{align}


\noindent Constraints (1.1) ensure that each school is assigned to exactly one start time. Constraints (1.2) ensure that every route for a given school arrives within the correct time window as defined by $\alpha$ and $\beta$. Constraints (1.3) set the minimum number of buses required to be greater than or equal to the number of routes in operation during any time period, as given in Proposition 1.  An analogous set of variables and constraints, detailed in Appendix A, implements combined morning and afternoon scheduling.

Together, constraints (1.1) and (1.2) imply that each route is assigned to at least one arrival time; in practice, the model will generally ensure that each route is assigned to exactly one arrival time in order to minimize $z$ via constraints (1.3). In the rare case where a route is assigned to multiple arrival times, any of these arrival times are feasible to achieve the minimum number of buses. While not strictly required, the following constraints can be added to ensure the model returns an optimal schedule with only one arrival time designated per bus:
\vspace{-0.9cm}

\begin{align}
&\sum_{m=1}^{M}x_{i, n}^{(m)}=1 & \forall n \in [N], i \in [\Gamma_{n}] \tag{1.6} 
\end{align}

The base model can easily be extended to incorporate practical scheduling considerations. For example, if certain schools are restricted from starting at certain times or are required to start at a certain time, constraints of the form $u_n^{(m)} = 0$ or $u_n^{(m)}= 1$, respectively, can be added as needed. Similar constraints for school ending time and bus arrival/departure can be added as well. If two schools are required to start at certain times relative to one another, constraints similar to (1.2) can be added as needed.

Every potential schedule affects the students, schools, and administration of a school district beyond the fixed cost of transportation associated with the schedule. As such, scheduling decisions cannot be based only on minimizing the number of buses required. We extend the base SBSP model to consider supplementary metrics, with a focus on equity, as additional criteria for determining school and route schedules.

\section{Equity and Start Time Changes}
We focus on the disutility associated with changing school start times as the equity determinant, and we pursue equity across schools. The formulations we present are easily adaptable to pursuing equity across other groups, including student demographic groups. The measure we examine is the \textit{center}, or maximum negative effect on any one group; minimizing this measure is known as the minimax method. A drawback of minimax is that it does not consider how equitably and efficiently effects are distributed among groups which do not experience the maximum negative effect (Luss, \citeyear{luss_lexminimax}). The lexicographic approach extends the minimax concept by sequentially minimizing the largest negative effect, followed by the second-largest, followed by the third-largest, and so on, leading to improved equitable solutions (Wang et al., \citeyear{wang}). Other measures of inequity reviewed in Marsh and Schilling (\citeyear{marsh_equity}) can be implemented within the model as required.

In addition, we use the concept of the ``price of fairness" as proposed by Bertsimas et al.\ (\citeyear{PoF}) to assess solutions. For a given schedule, let $f(\cdot)$ represent the value of some efficiency metric associated with the schedule. For example, $f(\cdot)$ could represent the number of buses required. For a given minimization problem, if OPT represents the most efficient scheduling solution in regard to $f(\cdot)$ and FAIR represents an equitable solution, then we define the price of fairness (PoF) as:
\vspace{-0.5cm}

$$\text{PoF(OPT, FAIR)} = \frac{f(\text{FAIR}) - f(\text{OPT})}{f(\text{OPT})} \ \ .$$

\noindent Bertsimas et al.\ (\citeyear{PoF}) consider the price of fairness to be generally unavoidable but potentially negligible based on the specific characteristics of the problem and utility set. We return to the idea of the price of fairness when analyzing the trade-offs between equity and efficiency in our numerical studies.

\subsection{Achieving Minimax Equity in Start Times}
For a school $n \in [N]$ and start time $m \in [M]$, let the parameter $c_{n}^{(m)}$ denote the disutility associated with changing school $n$'s current start time to $m$. The disutility can be positive or negative - a negative disutility implies that changing a school's start time benefits the school and its students. We define the term generally such that it is possible that there be a disutility associated with maintaining a current start time (e.g.\ the current start time of school $n$ is $m$, and $c_n^{(m)}$ is positive).

In our motivating work, we assume that the school district prefers schedules with lower absolute deviation in school start times from the current schedule. A proposed schedule is more likely to be approved by the school board and to be well-received in the community if the absolute changes in start times are distributed evenly across schools. In such a case, the disutility associated with changing school $n$'s start time from its current time $\widehat{m}$ to a new time $m$ can be represented as $c_n^{(m)} = |\widehat{m} - m|$. An efficient solution would minimize the sum of these changes, as in F{\"u}genschuh and Martin (\citeyear{fugenschuh2006}), but may lead to certain schools having significantly larger start time changes than others. 

We first pursue equity by seeking to minimize the maximum disutility associated with changing start times across all schools. This minimax approach, which promotes equity and efficiency concurrently, has the dual benefits of lending itself to straightforward ILP modeling and being easy to interpret from a policymaking perspective. We introduce a new variable $\pi_n$ to denote the disutility associated with school $n$'s start time change. Let the variable $\pi_{max}$ denote the maximum disutility across all schools and let $\varphi$ be a nonnegative scaling parameter between the two objectives. When $\varphi = 0$, the problem reduces to the basic SBSP from Section 3, whereas higher values of $\varphi$ correspond to a greater importance on $\pi_{max}$ when selecting a schedule. The SBSP with minimax equity across schools in regard to start time changes can be formulated as follows:
\vspace{-1.2cm}

\begin{align}
&min~~ z + \varphi\pi_{max}\nonumber \tag{ILP 2a}\\
s.t.~~&\pi_n = \sum_{m=\beta + 1}^M c_n^{(m)}u_n^{(m)}&\forall n \in [N] \nonumber \tag{2.1}\\
&\pi_{max} \geq \pi_n &\forall n \in [N] \nonumber \tag{2.2}\\
& (1.1)-(1.13) \textstyle \nonumber
\end{align}
\vspace{-0.9cm}

\noindent We simultaneously minimize the number of buses required and the maximum disutility. The relative importance of the two objectives is determined by the chosen value of $\varphi$. Constraints (2.1) calculate each school's associated disutility based on its new start time, and constraints (2.2) define $\pi_{max}$ as the maximum across all $\pi_n$. Defining the disutilities a priori in conjunction with the use of a time-indexed formulation allows for greater model flexibility. For example, a school district may determine start time change disutilities differently for each school, possibly leading to negative values or using complex nonlinear formulae. However the individual disutilities are defined, the model remains linear and the number of constraints needed to determine the disutility for each school do not increase. 

An alternate approach to the minimax problem is to minimize the maximum disutility given a number of available buses $\bar z$, an application of the $\epsilon$-constraint method for multi-objective optimization (Chankong and Haimes, \citeyear{chankong}):
\vspace{-0.9cm}

\begin{align}
&min~~ \pi_{max}\nonumber \tag{ILP 2b}\\
s.t.~~&z = \bar z &\tag{2.3}\\
& (1.1)-(1.13), (2.1), (2.2) \textstyle \nonumber
\end{align}

\subsection{Lexicographic Minimax Improvement}

Minimax focuses on minimizing the highest disutility across all schools, but does not guarantee an equitable or efficient distribution of disutility beyond the schools associated with the largest disutility. As with the base SBSP model, there may be many optimal solutions to the minimax problem. We propose a lexicographic minimax strategy to improve the efficiency of minimax solutions while maintaining or improving equity.

The lexicographic minimax approach, in the context of school scheduling, seeks to first minimize the highest disutility across schools, then the second highest disutility, etc. until finally minimizing the lowest ($N$th highest) disutility. Suppose we have two solutions, $A_1$ and $A_2$, to the minimax problem each requiring the same number of buses. Solution $A_1$ is lexicographically smaller than $A_2$ if and only if the $j$th highest disutility across schools for $A_1$ is less than or equal to that of $A_2$ for all $j \in [N]$ (Wang et al., \citeyear{wang}). Our strategy is to first solve either (ILP 2a) or (ILP 2b) to find the optimal number of buses $z^*$ and maximum disutility $\pi_{max}^*$. We then lexicographically improve the solution without increasing the number of buses or maximum disutility. The lexicographic minimax solution will be an optimal solution to the original minimax problem, equitable, and Pareto optimal across schools (Luss, \citeyear{luss_lexminimax}).

Consider the disutilities associated with each school ranked in descending order. For all $j \in \{1, \ldots, N\}$, let the variable $\psi_j$ represent the $j$th highest ranked disutility. Define the binary variable $h_{j, n}$ for all $j, n \in N$ to equal 1 if $\pi_n \leq \psi_j$ and 0 otherwise. Let the parameter $\widetilde{c}$ represent the largest difference between any two disutilities $c_{n_1}^{(m_1)}, c_{n_2}^{(m_2)}$. We formulate the lexicographic minimax problem as follows:
\vspace{-0.9cm}
 
\begin{align}
&lexmin~~ [\psi_1, \psi_2, \psi_3, ... , \psi_N] \nonumber \tag{ILP 3}\\
s.t.~~&z = z^* \nonumber \tag{3.1}\\
&\psi_1 = \pi_{max}^* \nonumber \tag{3.2}\\ \displaybreak 
&\sum_{n=1}^N h_{j, n} \geq N + 1 - j & \forall j \in [N] \tag{3.3} \\
& \widetilde{c} (1 - h_{j, n}) \geq \pi_n - \psi_j &\forall j,n \in [N] \tag{3.4} \\ 
& \psi_j \geq \psi_{j+1} & \forall j \in [N-1] \tag{3.5} \\
& h_{j, n} \in \{0, 1\} & \forall j, n \in [N] \tag{3.6} \\
& (1.1)-(1.13), (2.1) \textstyle \nonumber
\end{align}
\vspace{-0.9cm}

\noindent The objective minimizes each of the $j$th highest ranked disutilities in order, i.e., higher values of $j$ correspond to lower ranks and lowered priority. Constraints (3.1) and (3.2) define the number of buses and maximum disutility based on the solution to the minimax problem. Constraints (3.3) require that, for each $j \in [N]$, the appropriate number of school utilities are at or below the $j$th highest disutility. Constraints (3.4) require $h_{j, n}$ to equal 0 if $\psi_j > \pi_n$. Constraints (3.5) ensure the correct ordering of the $\psi_j$ variables.

One approach to solving this model directly is to scale each objective by a small positive value such that the priority order is preserved, i.e. $min \sum_{j = 1}^{N} \epsilon^n \psi_{j}$ where $\epsilon$ is sufficiently small. Attempting this approach with optimization software can lead to numerical issues. While Gurobi (\citeyear{gurobiguidelines}) recommends that the ratio between the largest and smallest coefficients in a model should not exceed $10^8$ or $10^9$, this guideline is violated even in a small instance with only ten schools and $\epsilon = 0.1$. 

Klotz and Newman (\citeyear{klotz}) recognize that large ratios -- by many orders of magnitude -- between MIP model coefficients may lead to round-off errors and inefficient solutions. As expected, Gurobi produces solutions inconsistent with the model when we test the approach on randomized instances with 15 schools. Specifically, (ILP 3) returns solutions with suboptimal $\psi_j$ for larger values of $j$. As an illustrative example, Gurobi solved a test instance with (ILP 3) to optimality (i.e., with an optimality gap of 0\%) in 76 seconds, producing the following solution:

\vspace{-0.4cm}
$$(\psi_{1}, \psi_{2}, \ldots, \psi_{15}) = (13, 11, 7, 7, 6, 6, 4, 4, 1, 1, 1, 1, 1, 1, 1)$$

\noindent However, lexicographically smaller solutions to the same test instance exist for which $\psi_{j} = 0$ for at least one $j$. Because the exponentially decreasing objective coefficients render the values of $\psi_j$ for larger $j$ negligible relative to the values of $\psi_j$ for smaller $j$, solvers are unable to distinguish between solutions which differ only on the $\psi_j$ for larger $j$. The issue of large ratios between model coefficients often requires attention when working with MIP models (see Talebian and Zou, (\citeyear{talebian}) and \c{C}a\u{g}lar and G{\"u}rel, (\citeyear{caglar}) for examples).

The solution to the lexicographic problem can instead be obtained by determining each $\psi_j$ iteratively:

\begin{algorithm}
	\SetKwInOut{Input}{Input}
	\SetKwInOut{Output}{Output}

	\Input{$z^*, \pi_{max}^*$ from (ILP 2), all relevant model parameters}
	\Output{Lexicographically optimal schedule}
	initialize $j^* = 2$;
	
	\While{$j^* \leq N$} 
	{
		solve (ILP 4) and return optimal objective value $\psi_{j^*}^*$;
		
		add constraint $\psi_{j^*} = \psi_{j^*}^*$ to (ILP 4);
		
		$j^* = j^* + 1$;
	}
	{
		return (ILP 4) solution for $j^* = N$
	}
	\caption{Iterative lexicographic schedule improvement}
\end{algorithm}


\noindent Without the constraints of the form $\psi_{j^*} = \psi_{j^*}^*$, (ILP 4) is as follows:
\begin{align}
\allowdisplaybreaks
&min~~ \psi_{j^*} \nonumber \tag{ILP 4}\\
s.t.~~&\sum_{n=1}^N h_{j, n} \geq N + 1 - j &\forall j \in [j^*]  \tag{4.1} \\
& \widetilde{c} (1 - h_{j, n}) \geq \pi_n - \psi_j &\forall j \in [j^*], n \in [N] \tag{4.2} \\ 
& h_{j, n} \in \{0, 1\} & \forall j \in [j^*], n \in [N] \tag{4.3} \\
& (1.1)-(1.13), (2.1), (3.1), (3.2) \textstyle \nonumber
\end{align}
\vspace{-0.9cm}

\noindent While the size of (ILP 4) increases with each iteration, an additional $\psi_j$ variable is fixed after each iteration. Thus, even the final iteration is easier to solve than (ILP 3).

It should be noted that some optimization software packages (Gurobi, CPLEX) provide functionality for optimizing a single model with lexicographically ordered objectives by assigning a priority level to each objective, while others (Xpress Solver, GLPK, LP\_SOLVE) do not. In the numerical study, we compare the performance of Algorithm 1 against Gurobi's built-in lexicographic optimization functionality, which uses a similar iterative process. 

\section{Numerical Studies}

In this section, we analyze our models using randomized instances based on real data from an urban public school district serving grades kindergarten through 8th. We compare the effectiveness of objectives in regard to start time change equity on instances of two sizes. Next, we assess the computational efficiency of implementing our models. Finally, we examine the managerial implications of objective trade-offs in the context of start time change equity. All models in this section are solved using Gurobi Optimizer 8.0.1 implemented in Python 3.6.6. 

\subsection{Scenarios I-III: 15 schools, 59 morning and afternoon routes}
We first construct and analyze instances with 15 schools serviced by 59 morning routes and 59 afternoon routes, the size of our partner school district's primary bus network. The district is among the largest 15\% of public school districts in the United States by population and number of schools (Gray et al., \citeyear{districtstats}). 

\begin{table}[!b]
	\small
	\centering 
	\begin{tabular}{ll}
		\hline \hline \\
		\textbf{Schools} & 15 \\
		&  \\
		\textbf{School length} & Randomly assigned 6 hr. 30 min., 6 hr. 45 min., or 7 hr. 30 min. \\
		&  \\
		\textbf{Total routes} & 59 morning, 59 afternoon \\
		&  \\
		\textbf{Routes per school} & Randomly distributed; equal number of morning and afternoon routes per school \\
		&  \\
		\textbf{Route length} & Randomly assigned between 10 and 50 min. plus 15 min. constant transition time \\
		&  \\
		\textbf{School start time window} & 8:15 AM to 9:20 AM \\
		&  \\
		\textbf{\begin{tabular}[c]{@{}l@{}}Morning route arrival\\ time window\end{tabular}} & Between 10 and 20 min. prior to school start; no arrivals prior to 7:55 AM \\
		&  \\
		\textbf{\begin{tabular}[c]{@{}l@{}}Afternoon route \\ departure time window\end{tabular}} & Between 5 and 20 min. after school start \\
		&  \\ \hline \hline
	\end{tabular}
\caption{Summary of parameters for randomized instances}
\end{table}

\subsubsection{Background and Parameters}
We construct ten randomized instances as detailed in Table 1. For modeling purposes, time is discretized by five minute intervals, and we discretize route lengths by rounding up to the next multiple of five minutes. As such, schools must start and end, and routes must arrive and depart on the five-minute marks relative to the hour (8:15 AM, 8:20 AM, 8:25 AM, etc.). In practice, all route transitions within the district occur in less than 15 minutes due to the geographic compactness of the district; the constant transition time provides an additional buffer in case of delays, increasing the reliability of the transportation system. As discussed earlier, we consider the disutility associated with changing a school's start time to be proportional to the absolute change in start time; i.e., the disutility associated with changing school $n$'s start time from $\widehat{m}$ to a new time $m$ is $c_n^{(m)} = |\widehat{m} - m|$.

For each of the ten instances, we consider three scenarios. In \textbf{Scenario I}, all schools currently start at 8:45 AM. In \textbf{Scenario II}, all schools currently start at 9:20 AM. In \textbf{Scenario III}, the current start times of the 15 schools are evenly split between 8:15 AM, 8:45 AM, and 9:20 AM. 

\subsubsection{Comparing Equity-Centered and Efficiency-Centered Models}
We initially apply the base model (including afternoon scheduling) to each instance, only minimizing the number of buses required, and observe the associated absolute changes in start time for each school. The optimal number of buses $z^*$ is 23 for two of the instances and 24 for the remaining eight instances. As expected, the start and end times of the schools vary greatly in the optimal solution to facilitate reusing buses. Figure 3 shows the distribution of start time changes for each scenario across all ten instances. 
\vspace{0.1cm}

\begin{figure}[h]
	\includegraphics{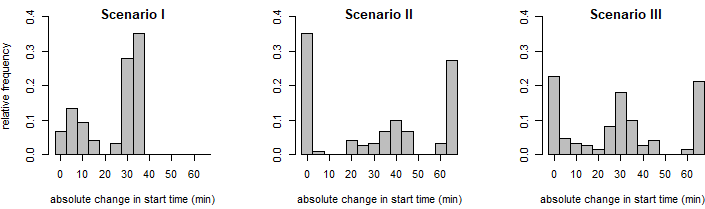} 
	\centering
	\caption{Start time change distribution, (ILP 1)}
\end{figure}
\vspace{-0.1cm}

\noindent These results illustrate the motivation for our work: focusing solely on reducing transportation costs in the modeling process leads to inefficient and inequitable distributions of start time changes. 

For every instance and scenario, we then set the number of available buses equal to $z^*$, as in (ILP 2b), and apply both our basic minimax and lexicographic minimax approaches. For comparison, we also optimize a ``Minimum Sum" objective, in which we set the number of available buses equal to $z^*$, then minimize the total disutility across all schools associated with changing or maintaining start times. It should be noted that minimizing the sum of disutilities is strictly an efficiency-based objective, unlike the minimax and lexicographic minimax objectives. We solve the problem with respect to each objective ten more times, setting the number of available buses equal to $z^*+ 1, z^* + 2, \ldots , z^* + 10$. Table 2 summarizes the average and standard deviation of start time changes within each instance for each scenario for each model. 
\vspace{0.4cm}

\begin{table}[!h]
\centering
\small
\renewcommand{\arraystretch}{1.5}
	\begin{tabular}{l|l|l||l|l||l|l||l|l|}
		\cline{2-9}
		& \multicolumn{2}{c||}{\textbf{Base}} & \multicolumn{2}{c||}{\textbf{Minimax}} & \multicolumn{2}{c||}{\textbf{Lex. Minimax}} & \multicolumn{2}{c|}{\textbf{Minimum Sum}} \\ \hline
		\multicolumn{1}{|c|}{Scen.} & 
		\multicolumn{1}{c|}{\begin{tabular}[c]{@{}c@{}}Avg.\\ Change\end{tabular}} & \multicolumn{1}{c||}{\begin{tabular}[c]{@{}c@{}}Std.\\ Dev.\end{tabular}} & \multicolumn{1}{c|}{\begin{tabular}[c]{@{}c@{}}Avg.\\ Change\end{tabular}} & \multicolumn{1}{c||}{\begin{tabular}[c]{@{}c@{}}Std.\\ Dev.\end{tabular}} & \multicolumn{1}{c|}{\begin{tabular}[c]{@{}c@{}}Avg.\\ Change\end{tabular}} & \multicolumn{1}{c||}{\begin{tabular}[c]{@{}c@{}}Standard\\ Dev.\end{tabular}} & \multicolumn{1}{c|}{\begin{tabular}[c]{@{}c@{}}Avg.\\ Change\end{tabular}} & \multicolumn{1}{c|}{\begin{tabular}[c]{@{}c@{}}Std. \\ Dev.\end{tabular}} \\ \hline
		\multicolumn{1}{|l|}{\textbf{I}: all 8:45 AM} & 23.80 & 12.59 & 24.57 & 11.88 & 21.07 & 11.81 & 19.77 & 13.17 \\ \hline
		\multicolumn{1}{|l|}{\textbf{II}: all 9:20 AM} & 31.60 & 26.44 & 30.77 & 26.70 & 26.00 & 24.35 & 24.47 & 25.04 \\ \hline
		\multicolumn{1}{|l|}{\textbf{III}: split times} & 29.67 & 22.54 & 22.57 & 18.31 & 19.17 & 17.28 & 16.90 & 18.74 \\ \hline
	\end{tabular}
\caption{Start time change results (in minutes), $z^*$ buses}
\end{table}

\noindent Appendix B contains similar tables for each scenario up to $z^* +10$ available buses as well as additional start time change distribution figures. We discuss key observations here.

Given a value of $z^*$, introducing the lexicographic minimax or minimum sum objective leads to lower average changes in absolute start time compared to the the base model, as expected. For all three objectives, the average change and the spread of the changes tend to decrease as the number of available buses increases. As an example, Figure 4 shows the relative frequency of absolute start time changes for the lexicographic minimax objective applied to Scenario III for $z^*$, $z^* + 4$, and $z^* + 8$ buses.
	
\begin{figure}[!h]
	\centering
	\includegraphics{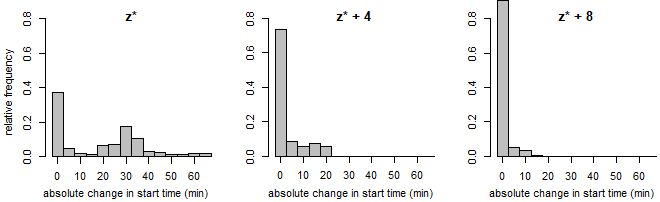}
	\caption{Start time change distribution, lexicographic minimax, Scenario III}
\end{figure}
\vspace{-0.3cm}

Compared to basic minimax, the lexicographic minimax solutions demonstrate significant improvements to both equity and efficiency in every scenario. The solutions produced by minimizing the sum tend to have a wider spread of absolute changes, with a higher maximum change but also more schools that do not change their start time at all.

\vspace{-0.3cm}
\begin{figure}[!h]
	\centering
	\includegraphics{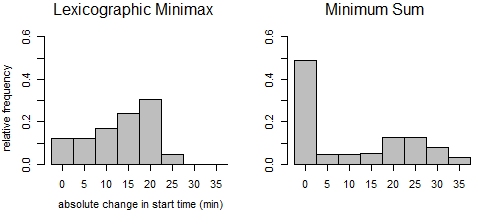}
	\caption{Start time change distribution, $z^* + 8$, Scenario I}
\end{figure}
\vspace{-0.3cm}

The lexicographic minimax approach should be used instead of basic minimax to improve equity in practical applications.The choice between the lexicographic minimax and minimum sum approaches is determined by whether equity or efficiency is prioritized in the planning process. As an example, Figure 5 shows the relative frequency of absolute start time changes for $z^* + 8$ buses in Scenario I. Compared to the lexicographic minimax solutions, the minimum sum solutions have significantly more schools with no start time change at the expense of some schools' start times changing by 30 or 35 minutes.

\subsubsection{Notes on Computational Efficiency}
The ten instances of the base model and the 330 minimax problems (ten instances with three scenarios each considering eleven possible values for the number of available buses) were solved to optimality without a time limit. Table 3 presents a summary of the time taken to solve these models.

\vspace{0.2cm}
\begin{table}[!h]
	\centering
	\begin{tabular}{l|l|l|l|l|l|l|}
		\cline{2-7}
		& Count & \begin{tabular}[c]{@{}l@{}}Mean CPU\\ Solve Time\end{tabular} & \begin{tabular}[c]{@{}l@{}}Median CPU\\ Solve Time\end{tabular} & \begin{tabular}[c]{@{}l@{}}Max CPU\\ Solve Time\end{tabular} & \begin{tabular}[c]{@{}l@{}}\% solved \\  in 1 min\end{tabular} & \begin{tabular}[c]{@{}l@{}}\% solved \\  in 5 min\end{tabular} \\ \hline
		\multicolumn{1}{|l|}{Base Model} & 10 & 62.53 sec & 18.52 sec & 3 min 39 sec & 70\% & 100\% \\ \hline
		\multicolumn{1}{|l|}{Minimax} & 330 & 85.03 sec & 10.85 sec & 137 min 16 sec & 91.5\% & 98.2\% \\ \hline
	\end{tabular}
	\caption{Summary of solve times for base and minimax model, Scenarios I-III}
\end{table}
\vspace{-0.2cm}

\noindent We observe that these problems are generally solvable to optimality in a short time using commercial optimization software for instances of comparable size. Every instance of the base model and nearly every minimax problem was solved in five minutes. In five of the six minimax cases which were not solved in five minutes, the optimal solution was found in under one minute (the last required 13 minutes and 45 seconds) with the vast majority of time spent on closing the optimality gap. 

We also observe that lexicographic minimax solutions can be computed efficiently using the iterative algorithm (Algorithm 1). Initial attempts to solve the lexicographic problem directly using (ILP 3) led to inconsistent solutions as expected. Thus, we used the iterative algorithm with a time limit of five minutes per iteration. To improve efficiency and to ensure a feasible solution was found within the time limit, we used the solution of the basic minimax problem as a warm start for the first iteration ($j^* = 2$), and we used the solution of the preceding iteration as a warm start for subsequent iterations. Of the 4,620 total iterations required to lexicographically improve all 330 basic minimax solutions, only 766 reached the time limit. 

The performance of Algorithm 1 is comparable to that of Gurobi's built-in lexicographic multiobjective functionality. Using the same five minute iteration time limit, 33 lexicographic improvement problems were re-solved using the built in functionality. Of these 33 problems solved using the built-in functionality, 7 produced better solutions than Algorithm 1, 6 produced worse solutions than Algorithm 1, and 20 produced identical solutions. On average, Algorithm 1 produced solutions 9.6\% faster than the built-in functionality. As such, Algorithm 1 (with warm starts) can be used to produce quality solutions to the lexicographic improvement problem.

Overall, the equity-centered models we present lead to improved solutions when compared to the base model. The lexicographic minimax approach outperforms the basic minimax model while remaining computationally tractable, and the lexicographic minimax approach leads to more equitable solutions than minimizing the sum in absolute changes in start times. 

\subsection{Scenario IV: 30 schools, 100 morning and afternoon routes}
To explore how these numerical and computational results hold for larger instances, we analyze a case with 30 schools serviced by 100 morning routes and 100 afternoon routes. Additionally, a strength of our time-indexed model formulations is that nonlinear disutility functions can be implemented without requiring additional constraints; we consider such a function here.

\subsubsection{Background and Parameters}

We consider five randomized instances of a newly constructed \textbf{Scenario IV}. Each instance of Scenario IV has 30 schools with a total of 100 morning routes and 100 afternoon routes. As in Scenario I, all schools currently start at 8:45 AM and schools are permitted to start only at five minute marks relative to the hour between 8:15 AM and 9:20 AM. 

We introduce a new disutility function in this scenario. For any school $n$ and for times $\widehat{m}, m \in [M]$, the disutility associated with changing school $n$'s start time from 
$\widehat{m}$ to $m$ is:
\vspace{-0.5cm}

\[ c_n^{(m)} = \begin{cases} 
(\widehat{m} - m)^{1.25} & \text{if } m < \widehat{m}, \\
0 & \text{if } m = \widehat{m}, \\
m - \widehat{m} & \text{if } m > \widehat{m}.
\end{cases}
\]

\noindent Figure 6 illustrates the disutility associated with changing a school's current start time to each potential new start time in the morning window.

\vspace{0.1in}
\begin{figure}[!h]
	\includegraphics[width=11cm]{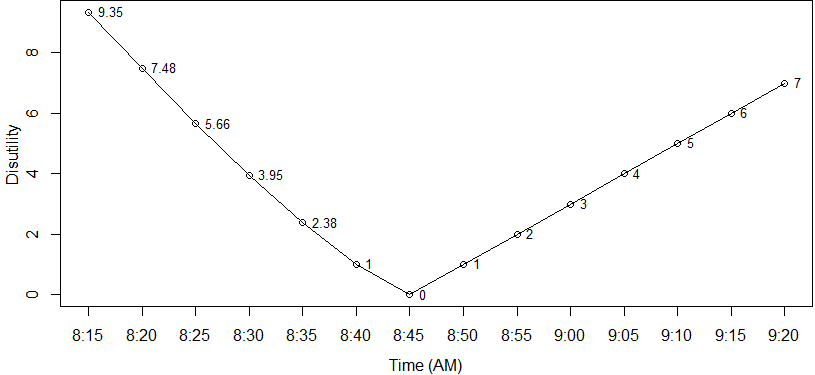} 
	\centering
	\caption{Disutility function values, Scenario IV}
\end{figure}
\vspace{-0.1in}

\noindent The remaining parameters are identical to those of Scenarios I-III (see Table 1).

\subsubsection{Comparing Equity-Centered and Efficiency-Centered Models}

For each of the five instances, we first solve the base model (ILP 1) to minimize the number of buses required. The minimum number of buses $z^*$ required in each instance to service the 100 morning and 100 afternoon routes is 41, 41, 42, 41, and 40, respectively. Figure 7  shows the distribution of start time change disutilities aggregated across all five instances.

\begin{figure}[h]
	\includegraphics[width=10cm]{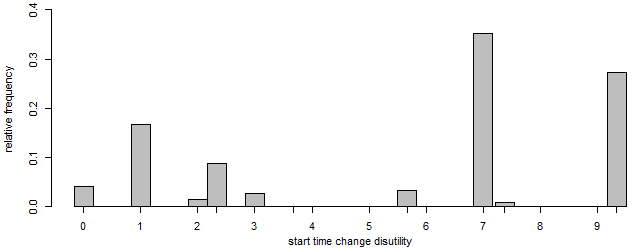} 
	\centering
	\caption{Scenario IV disutilities, (ILP 1)}
\end{figure} 
\vspace{-0.1in}

\noindent As in Scenarios I-III, only minimizing the number of buses required leads to inefficient and inequitable distributions of start time change disutilities. For each instance, we optimize in regard to the basic minimax, lexicographic minimax, and minimum sum objectives for $z^*$ and $z^* + 5$ available buses. Figure 8 and Figure 9 show the distribution of start time change disutilities aggregated across all five instances for $z^*$ and $z^* + 5$ buses, respectively.

\begin{figure}[!hp]
	\includegraphics{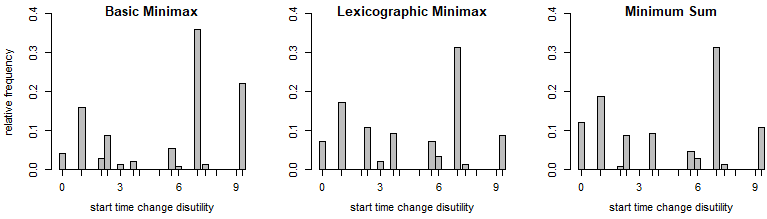} 
	\centering
	\caption{Scenario IV disutilities, $z^*$ buses}
\end{figure}
\vspace{0.1in}

\begin{figure}[!hp]
	\includegraphics{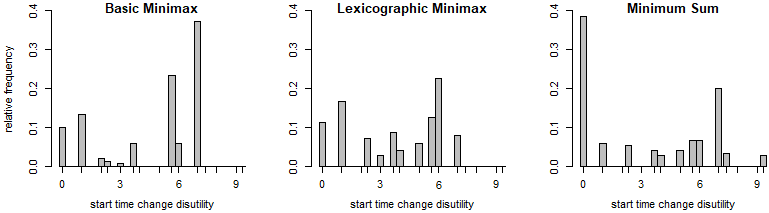} 
	\centering
	\caption{Scenario IV disutilities, $z^* + 5$ buses}
\end{figure}

For both $z^*$ and $z^* + 5$ available buses, the lexicographic minimax objective outperforms the basic minimax objective, just as in Scenarios I-III. We observe that for $z^*$ available buses, the disutility distributions for lexicographic minimax and minimum sum are similar, with the minimum sum objective function leading to slightly more schools experiencing disutilities at either end of the spectrum. This effect is stronger for $z^* + 5$ available buses, where minimum sum leads to more schools experiencing no start time change disutility at the cost of more schools experiencing the highest levels of start time change disutility as compared to lexicographic minimax. These observations reflect the trends observed in Scenarios I-III.

\subsubsection{Notes on Computational Efficiency}

For Scenario IV, the five instances of the base model and the ten minimax problems ($z^*$ and $z^* + 5$ for each instance) were solved to optimality. For more information on cases of this size, we also solved the minimax problem for $z^* + 1, z^* + 2, z^* + 3$ and $z^* + 4$. Table 4 presents a summary of the solve times.

\vspace{0.2cm}
\begin{table}[!h]
	\centering
	\begin{tabular}{l|l|l|l|l|l|l|}
		\cline{2-7}
		& Count & \begin{tabular}[c]{@{}l@{}}Mean CPU\\ Solve Time\end{tabular} & \begin{tabular}[c]{@{}l@{}}Median CPU\\ Solve Time\end{tabular} & \begin{tabular}[c]{@{}l@{}}Max CPU\\ Solve Time\end{tabular} & \begin{tabular}[c]{@{}l@{}}\% solved \\  in 5 min\end{tabular} & \begin{tabular}[c]{@{}l@{}}\% solved \\  in 15 min\end{tabular} \\ \hline
		\multicolumn{1}{|l|}{Base Model} & 5 & 22.73 sec & 11.79 sec & 59.85 sec & 100\% & 100\% \\ \hline
		\multicolumn{1}{|l|}{Minimax} & 30 & 252.29 sec & 110.88 sec & 32 min & 76.7\% & 96.7\% \\ \hline
	\end{tabular}
	\caption{Summary of solve times for base and minimax model}
\end{table}
\vspace{-0.2cm}

\noindent We observe that the base models of this size (30 schools, 100 morning and afternoon routes) require similar computational effort to solve as the smaller instances. While the minimax problems of this size remain tractable using commercial software, they require significantly longer computational time than the minimax problems of smaller size (a median of 110.88 sec vs. 10.85 sec). For instances which take longer to solve, much of the computational time is spent on closing the optimality gap, suggesting that future work should prioritize developing stronger lower bounds for the minimax problem.

We again use the iterative algorithm with warm starts to solve the lexicographic minimax problem for each instance with $z^*$ and $z^* + 5$ buses. Based on the computational results of the minimax problems, we use an iteration time limit of 15 minutes. Of the 290 total iterations required to lexicographically improve the 290 basic minimax, only 46 reached the time limit.

\subsection{Objective Trade-offs and Managerial Implications}
Previously, we fixed the number of buses $z$ available and minimized inequity via $\pi_{max}$. We now examine the managerial implications of pursuing various levels of equity in different scenarios by fixing $\pi_{max}$ and minimizing $z$. In this analysis, we include additional practical constraints reflecting a restriction faced by our partner school district: the maximum number of distinct school start times in the district cannot exceed three. Let the binary variable $k^{(m)}$ equal 1 if schools are permitted to start at time $m$ and 0 otherwise. We model the additional restriction as follows:
\vspace{-0.9cm}

\begin{align}
& \sum_{m=\beta+1}^{M}k^{(m)} \leq 3 & \tag{1.14} \\
& u_{n}^{(m)} \leq k^{(m)} &\forall n \in [N], m \in \{\beta+1, \ldots, M\} \tag{1.15}\\ 
&k^{(m)} \in \{0,1\} &\forall m \in \{\beta+1, \ldots, M\} \tag{1.16}
\end{align}

\noindent Constraints (1.14) limit the total number of distinct start times to three. Constraints (1.15) ensure that schools must start at one of the distinct start times as defined by the $k^{(m)}$ variables.

For an instance, let OPT represent the solution which minimizes the number of buses required with no restriction on the changes in start time. Let FAIR$_{\tau}$ represent the solution which minimizes the number of buses required subject to the maximum absolute change in start time across all schools being bounded above by $\tau$ minutes. For any solution, let $f(\cdot)$ represent the required number of buses associated with the solution. Note that $f(\text{FAIR}_0)$ is the minimum number of buses required to service the current school bell time schedule. For every instance in every scenario, we solve for FAIR$_0$, FAIR$_5$, FAIR$_{10}$, $\ldots$ \  and OPT. Figure 10 displays the average results for Scenario I across all ten instances.

\begin{figure}[!h]
	\includegraphics[width=9cm]{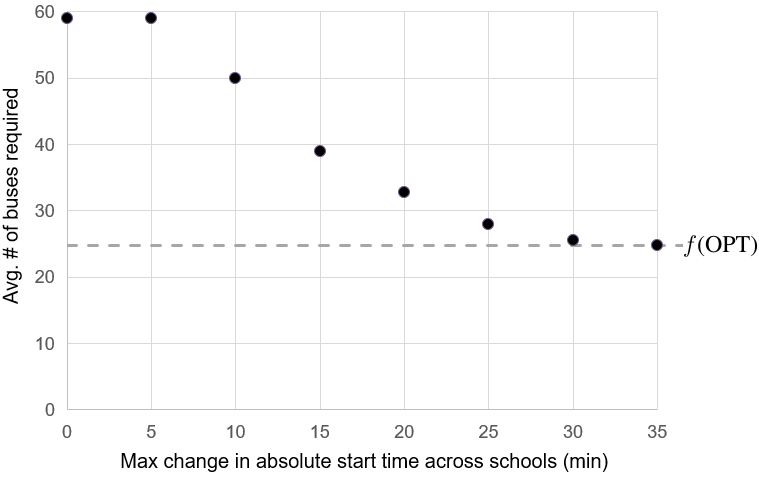} 
	\centering
	\caption{$f(\text{FAIR}_\tau), \tau = \{0, 5, \ldots, 35\}$, Scenario I}
\end{figure}

We first observe that the price of fairness PoF(OPT, FAIR$_\tau$) is large for smaller values of $\tau$. The key characteristic of solutions requiring fewer of buses is that school start \textit{and} end times are spread out. However, because all schools currently start at the same time in Scenario I, smaller values of $\tau$ restrict the spread of start times. At the same time, these results are favorable from an alternate managerial perspective. A district can equitably save a significant number of school buses with a maximum absolute change in school start time of as little as 15 minutes, translating to a potential annual savings of over \$750,000 (based on the bus rental rates of our partner district) with a slight change in schedule which may face lower resistance from the district's community. The results for Scenario II, displayed in Figure 11, exhibit some differences.

\begin{figure}[!h]
	\includegraphics[width=9cm]{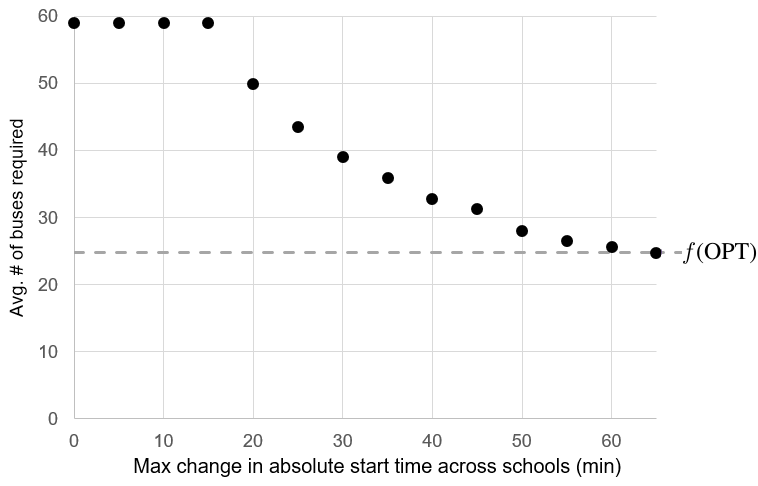} 
	\centering
	\caption{$f(\text{FAIR}_\tau), \tau = \{0, 5, \ldots, 65\}$, Scenario II}
\end{figure}

Like Scenario I, in Scenario II the price of fairness PoF(OPT, FAIR$_\tau$) is large for smaller values of $\tau$. However, the price of fairness does not decrease until $\tau =$ 20 minutes. Consider that Scenario I, where all schools current start at 8:45, allows a spread of start times between 8:30 AM and 9:00 AM. The same $\tau =$ 15 minutes in in Scenario II, a $\tau$ value of 15 minutes, necessitates every school to begin between 9:05 AM and 9:20 AM because all schools currently start at 9:20 AM, the end of the allowable start time window. As such, it requires sacrificing equity and efficiency in start time changes to achieve a significant reduction in the number of buses required. Figure 12 displays the results for Scenario III.

\begin{figure}[!h]
	\includegraphics[width=9cm]{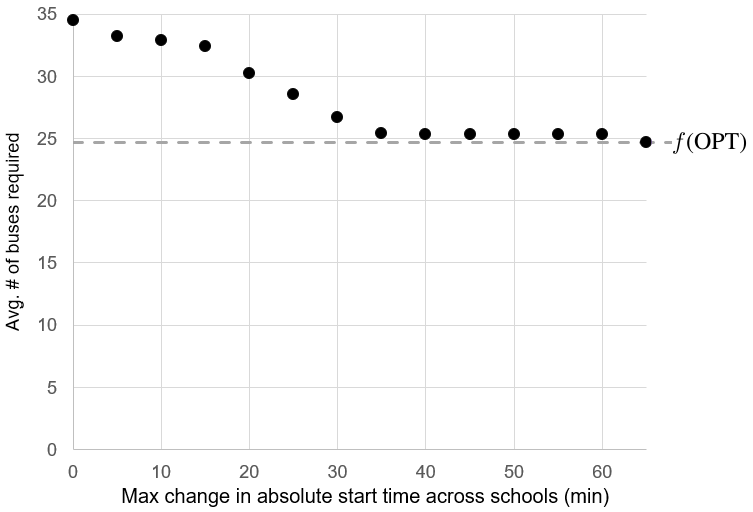} 
	\centering
	\caption{$f(\text{FAIR}_\tau), \tau = \{0, 5, \ldots, 65\}$, Scenario III}
\end{figure}

The key observation is that the price of fairness PoF(OPT, FAIR$_\tau$) is smaller in this scenario for all values of $\tau$ as the number of buses required with the current schedule is much lower on account of the current start times being spread across the interval. However, this implies that it is difficult to achieve a significant reduction in the number of buses while maintaining an equitable and efficient distribution in school start time changes when school start and end times are already spread out.

\section{Concluding Remarks}

In this paper, we explicitly consider equity as a factor when determining optimal school bell times and route schedules for public school transportation systems. Motivated in part through a collaboration with a moderately-sized public school district, we also seek to extend the focus on equity in public sector operations research to the school bus scheduling problem. We present a minimax formulation for equitably reducing the disutilities associated with changing school start times. We then propose a lexicographic minimax strategy with an iterative algorithm for improving minimax solutions. 

We perform numerical studies on randomized instances based on actual data from a public school transportation system. Quantifying disutility as the absolute change in a school's start time, we compare the efficiency and equity of schedules produced by the base model, by minimax, by lexicographic minimax, and by minimizing the total disutility. We find that the lexicographic model significantly outperforms the base and minimax models, and produces more equitable solutions than simply minimizing total disutility.  Our numerical studies demonstrate the practical applicability of our models to improving public school transportation systems using optimization software.

A limitation of our approach is that we do not consider route-dependent or location-dependent bus transition times between routes. Instead, we use constant transition times due to the geographically compact nature of our collaborating school district, so the time-indexed formulation proposed by Zeng et al.\ (\citeyear{zeng_sbsp}) serves as an appropriate foundation for our analysis. However, the general strategies we use to address issues of equity can be readily adapted to formulations which use exact transition times, presenting a potential avenue for future work. Future research should also focus on developing provable bounds and heuristic algorithms for our models, as relying on commercial optimization software may be impractical for instances significantly larger than the ones analyzed in this paper.

\section*{Acknowledgments} 
This project is supported by the National Science Foundation (CMMI-1727744) and the Northwestern University Transportation Center (Greenbriar Equity Group Fellowship).

  
\bibliography{sbsp_bib}

\begin{thebibliography}{46}
\providecommand{\natexlab}[1]{#1}
\providecommand{\url}[1]{\texttt{#1}}
\expandafter\ifx\csname urlstyle\endcsname\relax
  \providecommand{\doi}[1]{doi: #1}\else
  \providecommand{\doi}{doi: \begingroup \urlstyle{rm}\Url}\fi

\bibitem[Bertsimas et~al.(2011)Bertsimas, Farias, and Trichakis]{PoF}
Dimitris Bertsimas, Vivek~F. Farias, and Nikolaos Trichakis.
\newblock The price of fairness.
\newblock \emph{Operations Research}, 59\penalty0 (1):\penalty0 17--31, 2011.

\bibitem[Bertsimas et~al.(2019)Bertsimas, Delarue, and Martin]{mit_sbsp}
Dimitris Bertsimas, Arthur Delarue, and Sebastien Martin.
\newblock Optimizing schools' start time and bus routes.
\newblock \emph{Proceedings of the National Academy of Sciences}, 116\penalty0
  (13):\penalty0 5943--5948, 2019.

\bibitem[B{\"o}gl et~al.(2015)B{\"o}gl, Doerner, and Parragh]{bogl}
Michael B{\"o}gl, Karl~F. Doerner, and Sophie~N. Parragh.
\newblock The school bus routing and scheduling problem with transfers.
\newblock \emph{Networks}, 65\penalty0 (2):\penalty0 180--203, 2015.

\bibitem[Bookbinder and Edwards(1990)]{bookbinder_SBSP}
James~H. Bookbinder and Steven~H. Edwards.
\newblock School-bus routing for program scheduling.
\newblock \emph{Computers \& Operations Research}, 17\penalty0 (1):\penalty0
  79--94, 1990.

\bibitem[Bouzarth et~al.(2018)Bouzarth, Forrester, Hutson, and
  Reddoch]{bouzarth_districting}
Elizabeth~L. Bouzarth, Richard Forrester, Kevin~R. Hutson, and Lattie Reddoch.
\newblock Assigning students to schools to minimize both transportation costs
  and socioeconomic variation between schools.
\newblock \emph{Socio-Economic Planning Sciences}, 64:\penalty0 1--8, 2018.

\bibitem[Bowerman et~al.(1995)Bowerman, Hall, and Calamai]{bowerman}
Robert Bowerman, Brent Hall, and Paul Calamai.
\newblock A multi-objective optimization approach to urban school bus routing:
  Formulation and solution method.
\newblock \emph{Transportation Research Part A}, 29\penalty0 (2):\penalty0
  107--123, 1995.

\bibitem[Campbell et~al.(2008)Campbell, Vandenbussche, and
  Hermann]{campbell_equity}
Ann~M. Campbell, Dieter Vandenbussche, and William Hermann.
\newblock Routing for relief efforts.
\newblock \emph{Transportation Science}, 42\penalty0 (2):\penalty0 127--145,
  2008.

\bibitem[\c{C}a\u{g}lar and G{\"u}rel(2017)]{caglar}
Musa \c{C}a\u{g}lar and Sina G{\"u}rel.
\newblock Public {R \& D} project portfolio selection problem with
  cancellations.
\newblock \emph{{OR} Spectrum}, 39\penalty0 (3):\penalty0 659--687, 2017.

\bibitem[Chankong and Haimes(1983)]{chankong}
Vira Chankong and Yacov~Y. Haimes.
\newblock \emph{Multiobjective decision making: {Theory} and methodology}.
\newblock North-Holland, New York, NY, 1983.

\bibitem[Chen et~al.(2015)Chen, Kong, Dang, Hou, and Ye]{chen}
Xiaopan Chen, Yunfeng Kong, Lanxue Dang, Yane Hou, and Xinyue Ye.
\newblock Exact and metaheuristic approaches for a bi-objective school bus
  scheduling problem.
\newblock \emph{{PLOS ONE}}, 10\penalty0 (7), 2015.

\bibitem[Delgado~Serna and Pacheco~Bonrostro(2001)]{delgado_SBRPequity}
Cristina~R. Delgado~Serna and Joaqu{\'i}n Pacheco~Bonrostro.
\newblock Minimax vehicle routing problems: Application to school transport in
  the province of {Burgos}.
\newblock In Stefan Vo\ss{} and Joachim~R. Daduna, editors,
  \emph{Computer-Aided Scheduling of Public Transport}, pages 297--317, 2001.

\bibitem[Desrosiers et~al.(1986)Desrosiers, Ferland, Rousseau, Lapalme, and
  Chapleau]{transcol}
Jacques Desrosiers, Jacques-A Ferland, Jean-Marc Rousseau, Guy Lapalme, and Luc
  Chapleau.
\newblock {TRANSCOL}: {A} multi-period school bus routing and scheduling
  system.
\newblock \emph{TIMS Studies in the Management Sciences}, 22:\penalty0 47--71,
  1986.

\bibitem[Dunietz et~al.(2017)Dunietz, Matos-Moreno, Singer, Davis, O'Brien, and
  Chervin]{dunietz_start_times}
Galit~Levi Dunietz, Amilcar Matos-Moreno, Dianne~C. Singer, Matthew~M. Davis,
  Louise~M. O'Brien, and Ronald~D. Chervin.
\newblock Later school start times: What informs parent support or opposition?
\newblock \emph{Journal of Clinical Sleep Medicine}, 13\penalty0 (7):\penalty0
  889--897, 2017.

\bibitem[F{\"u}genschuh(2009)]{fugenschuh}
Armin F{\"u}genschuh.
\newblock Solving a school bus scheduling problem with integer programming.
\newblock \emph{European Journal of Operational Research}, 193\penalty0
  (3):\penalty0 867--884, 2009.

\bibitem[F{\"u}genschuh(2011)]{fugenschuh2011}
Armin F{\"u}genschuh.
\newblock A set partitioning reformulation of a school bus scheduling problem.
\newblock \emph{Journal of Scheduling}, 14\penalty0 (4):\penalty0 307--318,
  2011.

\bibitem[F{\"u}genschuh and Martin(2006)]{fugenschuh2006}
Armin F{\"u}genschuh and Alexander Martin.
\newblock A multicriteria approach for optimizing bus schedules and school
  starting times.
\newblock \emph{Annals of Operations Research}, 147\penalty0 (1):\penalty0
  199--216, 2006.

\bibitem[Gray et~al.(2013)Gray, Bitterman, and Goldring]{districtstats}
Lucinda Gray, Amy Bitterman, and Rebecca Goldring.
\newblock Characteristics of public school districts in the {United States}:
  Results from the 2011-12 schools and staffing survey ({NCES} 2013-311).
\newblock Technical report, National Center for Education Statistics, 2013.
\newblock URL \url{https://nces.ed.gov/pubs2013/2013311.pdf}.
\newblock Accessed 27 Nov.\ 2018.

\bibitem[Gur(2017)]{gurobiguidelines}
\emph{{Gurobi Guidelines for Numerical Issues}}.
\newblock Gurobi Optimization, 2017.
\newblock URL \url{http://files.gurobi.com/Numerics.pdf}.
\newblock Accessed 27 Nov.\ 2018.

\bibitem[Hashi et~al.(2016)Hashi, Hasan, and Zaman]{hashi}
Emrana~K. Hashi, M.R. Hasan, and M.S.U. Zaman.
\newblock {GIS} based heuristic solution of the vehicle routing problem to
  optimize the school bus routing and scheduling.
\newblock In \emph{2016 19th International Conference on Computer and
  Information Technology (ICCIT)}, pages 56--60. IEEE, 2016.

\bibitem[Johnson et~al.(2018)Johnson, Midgley, and Chichirau]{johnson}
Michael~P. Johnson, Gerald Midgley, and George Chichirau.
\newblock Emerging trends and new frontiers in community operational research.
\newblock \emph{European Journal of Operational Research}, 268\penalty0
  (3):\penalty0 1178--1191, 2018.

\bibitem[Kim et~al.(2012)Kim, Kim, and Park]{kim_kim_park}
Byung-In Kim, Seongbae Kim, and Junhyuk Park.
\newblock A school bus scheduling problem.
\newblock \emph{European Journal of Operational Research}, 218\penalty0
  (2):\penalty0 577--585, 2012.

\bibitem[Klotz and Newman(2013)]{klotz}
Ed~Klotz and Alexandra~M. Newman.
\newblock Practical guidelines for solving difficult mixed integer linear
  programs.
\newblock \emph{Surveys in Operations Research and Management Science},
  18\penalty0 (1--2):\penalty0 18--32, 2013.

\bibitem[Leachman et~al.(2017)Leachman, Masterson, and Figueroa]{punishing}
Michael Leachman, Kathleen Masterson, and Eric Figueroa.
\newblock A punishing decade for school funding.
\newblock Technical report, Center on Budget and Policy Priorities, 2017.
\newblock URL
  \url{https://www.cbpp.org/sites/default/files/atoms/files/11-29-17sfp.pdf}.
\newblock Accessed 27 Nov.\ 2018.

\bibitem[Leclerc et~al.(2012)Leclerc, McLay, and Mayorga]{leclerc}
Philip~D. Leclerc, Laura~A. McLay, and Maria~E. Mayorga.
\newblock Modeling equity for allocating public resources.
\newblock In Michael~P. Johnson, editor, \emph{Community-Based Operations
  Research}, pages 97--118, 2012.

\bibitem[Li and Fu(2002)]{li_SBRPequity}
Leon Li and Zhuo Fu.
\newblock The school bus routing problem: a case study.
\newblock \emph{Journal of the Operational Research Society}, 53\penalty0
  (5):\penalty0 552--558, 2002.

\bibitem[Litman(2002)]{litman_equity}
Todd Litman.
\newblock Evaluating {Transportation} {Equity}.
\newblock \emph{World Transport Policy \& Practice}, 8\penalty0 (2):\penalty0
  50--65, 2002.

\bibitem[Luss(1999)]{luss_lexminimax}
Hanan Luss.
\newblock On equitable resource allocation problems: A lexicographic minimax
  approach.
\newblock \emph{Operations Research}, 47\penalty0 (3):\penalty0 361--378, 1999.

\bibitem[Marsh and Schilling(1994)]{marsh_equity}
Michael~T. Marsh and David~A. Schilling.
\newblock Equity measurement in facility location analysis: {A} review and
  framework.
\newblock \emph{European Journal of Operational Research}, 74\penalty0
  (1):\penalty0 1--17, 1994.

\bibitem[Olariu(1991)]{olariu}
Stephan Olariu.
\newblock An optimal greedy heuristic to color interval graphs.
\newblock \emph{Information Processing Letters}, 37\penalty0 (1):\penalty0
  21--25, 1991.

\bibitem[Owens et~al.(2010)Owens, Belon, and Moss]{owens}
Judith~A. Owens, Katherine Belon, and Patricia Moss.
\newblock Impact of delaying school start time on adolescent sleep, mood, and
  behavior.
\newblock \emph{Archives of Pediatrics and Adolescent Medicine}, 164\penalty0
  (7):\penalty0 608--614, 2010.

\bibitem[Pacheco et~al.(2013)Pacheco, Caballero, Laguna, and
  Molina]{pacheco_SBRPequity}
Joaqu{\'i}n Pacheco, Rafael Caballero, Manuel Laguna, and Juli{\'a}n Molina.
\newblock Bi-objective bus routing: An application to school buses in rural
  areas.
\newblock \emph{Transportation Science}, 47\penalty0 (3):\penalty0 397--411,
  2013.

\bibitem[Park and Kim(2010)]{park_kim_review}
Junhyuk Park and Byung-In Kim.
\newblock The school bus routing problem: {A} review.
\newblock \emph{European Journal of Operational Research}, 202\penalty0
  (2):\penalty0 311--319, 2010.

\bibitem[Racker(2018)]{california}
Mini Racker.
\newblock California {Gov}.\ {Jerry} {Brown} rejects bill to prohibit schools
  from starting before 8:30 a.m.
\newblock \emph{Los Angeles Times}, 2018.
\newblock URL
  \url{http://www.latimes.com/politics/la-pol-ca-school-start-time-vetoed-20180921-story.html}.
\newblock Accessed 27 Nov.\ 2018.

\bibitem[Raff(1983)]{raff}
Samuel Raff.
\newblock Routing and scheduling of vehicles and crews: The state of the art.
\newblock \emph{Computers \& Operations Research}, 10\penalty0 (2):\penalty0
  63--211, 1983.

\bibitem[Rowe(2016)]{missouricuts}
Shelby Rowe.
\newblock School budget cuts mostly affect transportation.
\newblock \emph{Fulton Sun}, 2016.
\newblock URL
  \url{http://www.fultonsun.com/news/local/story/2016/sep/20/school-budget-cuts-mostly-affect-transportation/640941/}.
\newblock Accessed 27 Nov.\ 2018.

\bibitem[Savas(1978)]{savas_equity}
E.S. Savas.
\newblock Equity in providing public services.
\newblock \emph{Management Science}, 24\penalty0 (8):\penalty0 800--808, 1978.

\bibitem[Shafahi et~al.(2018)Shafahi, Wang, and Haghani]{shafahi}
Ali Shafahi, Zhongxiang Wang, and Ali Haghani.
\newblock {SpeedRoute}: Fast, efficient solutions for school bus routing
  problems.
\newblock \emph{Transportation Research Part B}, 117:\penalty0 473--493, 2018.

\bibitem[Spada et~al.(2005)Spada, Bierlaire, and Liebling]{spada}
Michela Spada, Michel Bierlaire, and Thomas~M. Liebling.
\newblock Decision-aiding methodology for the school bus routing and scheduling
  problem.
\newblock \emph{Transportation Science}, 39\penalty0 (4):\penalty0 477--490,
  2005.

\bibitem[Spears(2018)]{kentuckycuts}
Valarie~H. Spears.
\newblock ‘{How} in the world can we make it?’ {School} districts stunned
  by transportation cut.
\newblock \emph{Lexington Herald-Leader}, 2018.
\newblock URL
  \url{https://www.kentucky.com/news/local/education/article195214849.html}.
\newblock Accessed 27 Nov.\ 2018.

\bibitem[Swersey and Ballard(1984)]{swersey}
Arthur~J. Swersey and Wilson Ballard.
\newblock Scheduling school buses.
\newblock \emph{Management Science}, 30\penalty0 (7):\penalty0 844--853, 1984.

\bibitem[Talebian and Zou(2015)]{talebian}
Ahmadreza Talebian and Bo~Zou.
\newblock Integrated modeling of high performance passenger and freight train
  planning on shared-use corridors in the {US}.
\newblock \emph{Transportation Research Part B}, 82:\penalty0 114--140, 2015.

\bibitem[{Urban Institute Student Transportation Working
  Group}(2017)]{national_costs}
{Urban Institute Student Transportation Working Group}.
\newblock Student transportation and educational access.
\newblock Technical report, Urban Institute, 2017.
\newblock Accessed 27 Nov.\ 2018.

\bibitem[Wang et~al.(2004)Wang, Fang, and Hipel]{wang}
Lizhong Wang, Liping Fang, and Keith~W. Hipel.
\newblock Lexicographic minimax approach to fair water allocation problems.
\newblock In \emph{2004 IEEE International Conference on Systems, Man and
  Cybernetics}, volume~1, pages 1038--1043. IEEE, 2004.

\bibitem[Wolfson et~al.(2007)Wolfson, Spaulding, and Baroni]{wolfson}
Amy~R. Wolfson, Noah~L. Spaulding, and Elizabeth~M. Baroni.
\newblock Middle school start times: The importance of a good night's sleep for
  young adolescents.
\newblock \emph{Behavioral Sleep Medicine}, 5\penalty0 (3):\penalty0 194--209,
  2007.

\bibitem[Yan et~al.(2015)Yan, Hsaio, and Chen]{yan}
Shangyao Yan, Fei-Yen Hsaio, and Yi-Chun Chen.
\newblock Inter-school bus scheduling under stochastic travel times.
\newblock \emph{Networks and Spatial Economics}, 15\penalty0 (4):\penalty0
  1049--1074, 2015.

\bibitem[Zeng et~al.(2018)Zeng, Chopra, and Smilowitz]{zeng_sbsp}
Liwei Zeng, Sunil Chopra, and Karen Smilowitz.
\newblock A bounded formulation for the school bus scheduling problem.
\newblock 2018.
\newblock URL \url{https://arxiv.org/pdf/1803.09040.pdf}.
\newblock In review, accessed 27 Nov.\ 2018.

\end{thebibliography}
\bibliographystyle{plainnat}     


\appendix
\setcounter{figure}{0}
\setcounter{table}{0}

\makeatletter 
\renewcommand{\thefigure}{B.\@arabic\c@figure}
\makeatother	

\makeatletter 
\renewcommand{\thetable}{B.\@arabic\c@table}
\makeatother

\section*{Appendix A. \hspace{0.1cm} Afternoon Scheduling Constraints}
Let $\delta_{n}$ denote the length of the school day for school $n$. Let $\Theta_{n}$ represent the number of afternoon routes for school $n$. The $i^{th}$ afternoon route associated with school $n$ has a discretized length (duration) of $r_{i,n}^{PM}$ consecutive time periods, including a constant transition time added to each route. Let the parameter $\lambda$ denote the minimum number of time periods after the ending time of a school that any associated route may depart from that school, and let $\mu$ denote the maximum number of time periods after the ending time of a school that all associated routes must depart from that school. We denote the set of potential bus departure times as $P = \{p_{min}, \ldots, p_{max}\}$; to prevent infeasibility, it should hold that $p_{min} \leq M + \min\{\delta_{n}\}$ and $p_{max} \geq M + \lambda + \max\{\delta_{n}\}$, as well as $p_{max} \geq p_{min} + \lambda$. 

For the afternoon counterparts to the morning scheduling variables, define binary variables $y_{i, n}^{(m)}$ for all $n \in [N], i \in \Gamma_{n}, m \in [M]$ and $v_{n}^{(m)}$ for all $n \in [N], m \in \{p_{min}, \ldots, p_{max}-\lambda_{n}\}$ as such: $y_{i, n}^{(m)} = 1$ if route $i \in \Gamma_{n}$ departs at time $m$ and 0 otherwise; $v_{n}^{(m)} = 1$ if school $n$ ends at $m$ and 0 otherwise. The following constraints implement afternoon scheduling:

\begin{align}
& \sum_{m=p_{min}}^{p_{max}-\lambda}v_{n}^{(m)} = 1 &\forall n \in [N] \tag{1.7} \\ 
& v_{n}^{(m)} \leq \sum_{t=m+\lambda}^{\min \{m+\mu, p_{max}\}}y_{i,n}^{(t)} &\forall n \in [N], i \in [\Theta_{n}], m \in \{p_{min}, \ldots, p_{max}-\lambda\} \nonumber \tag{1.8}\\ 
&\sum_{n=1}^{N}\sum_{i=1}^{\Theta_{n}}\sum_{t=\max \{m-r_{i, n}^{PM}+1, p_{min}\}}^{m}y_{i,n}^{(t)}\leq z &\forall m \in P \tag{1.9}\\
& u_{n}^{(m)} = v_{n}^{(m+\delta_n)} & \forall n \in [N], m \in \{\max \{1, p_{min}-\delta_n\}, \ldots, \nonumber \\
&& \min \{M, p_{max} - \delta_n \} \}\tag{1.10} \\
&y_{i,n}^{(m)} \in \{0,1\} &\forall n \in [N], i \in [\Theta_{n}], m \in P \tag{1.11} \\
&v_{n}^{(m)} \in \{0,1\} &\forall n \in [N], m \in \{p_{min}, \ldots, p_{max}-\lambda\} \tag{1.12}
\end{align}

\vspace{0.3cm} 

\noindent Constraints (1.7) ensure that each school is assigned to exactly one start time. Constraints (1.8) ensure that every morning route for a given school arrives within the correct time window as defined by $\lambda$ and $\mu$. Constraints (1.9) set the minimum number of buses required to be greater than or equal to the number of routes in operation during any time period, as given in Proposition 1. Constraints (1.10) link the otherwise independent morning and afternoon variables using $\delta_n$.

The following optional constraints, analogous to (1.6), require each route to be assigned to exactly one departure time:
\vspace{-0.9cm}

\begin{align}
&\sum_{m=p_{min}}^{p_{max}}y_{i, n}^{(m)}=1 & \forall n \in [N], i \in [\Theta_{n}] \tag{1.13} 
\end{align}

\vspace{1cm}  
\section*{Appendix B. \hspace{0.1cm} Results of Start Time Change Analysis}
\vspace{0.5cm}
\begin{table}[!h]
	\centering
	\begin{tabular}{c|l|l||l|l||l|l|}
		\cline{2-7}
		\multicolumn{1}{l|}{} & \multicolumn{2}{c||}{\textbf{Minimax}} & \multicolumn{2}{c||}{\textbf{Lexicographic Minimax}} & \multicolumn{2}{c|}{\textbf{Minimum Sum}} \\ \hline
		\multicolumn{1}{|c|}{Buses} & \multicolumn{1}{c|}{\begin{tabular}[c]{@{}c@{}}Average\\ Change\end{tabular}} & \multicolumn{1}{c||}{\begin{tabular}[c]{@{}c@{}}Standard\\ Deviation\end{tabular}} & \multicolumn{1}{c|}{\begin{tabular}[c]{@{}c@{}}Average\\ Change\end{tabular}} & \multicolumn{1}{c||}{\begin{tabular}[c]{@{}c@{}}Standard\\ Deviation\end{tabular}} & \multicolumn{1}{c|}{\begin{tabular}[c]{@{}c@{}}Average\\ Change\end{tabular}} & \multicolumn{1}{c|}{\begin{tabular}[c]{@{}c@{}}Standard \\ Deviation\end{tabular}} \\ \hline
		\multicolumn{1}{|l|}{$z^*$} & 24.57 & 11.88 & 21.07 & 11.81 & 19.77 & 13.17 \\ \hline
		\multicolumn{1}{|l|}{$z^*+1$} & 23.60 & 10.94 & 19.70 & 10.56 & 16.97 & 13.91 \\ \hline
		\multicolumn{1}{|l|}{$z^*+2$} & 23.13 & 10.60 & 18.10 & 9.41 & 15.53 & 13.89 \\ \hline
		\multicolumn{1}{|l|}{$z^*+3$} & 22.67 & 9.46 & 17.27 & 8.60 & 14.60 & 13.83 \\ \hline
		\multicolumn{1}{|l|}{$z^*+4$} & 21.53 & 8.37 & 16.53 & 8.31 & 13.77 & 13.31 \\ \hline
		\multicolumn{1}{|l|}{$z^*+5$} & 20.83 & 7.07 & 15.60 & 8.16 & 12.87 & 12.81 \\ \hline
		\multicolumn{1}{|l|}{$z^*+6$} & 20.10 & 6.62 & 14.50 & 7.64 & 12.17 & 12.58 \\ \hline
		\multicolumn{1}{|l|}{$z^*+7$} & 18.90 & 6.01 & 14.13 & 7.21 & 11.40 & 12.47 \\ \hline
		\multicolumn{1}{|l|}{$z^*+8$} & 18.07 & 5.97 & 13.17 & 7.11 & 10.77 & 12.21 \\ \hline
		\multicolumn{1}{|l|}{$z^*+9$} & 17.53 & 4.95 & 12.97 & 6.43 & 10.10 & 11.96 \\ \hline
		\multicolumn{1}{|l|}{$z^*+10$} & 17.17 & 4.51 & 12.40 & 6.14 & 9.50 & 11.58 \\ \hline
	\end{tabular}
\caption{Average and standard deviation of absolute start time changes within each instance (in minutes), Scenario I}
\end{table}

\begin{table}[!h]
	\centering
	\begin{tabular}{c|l|l||l|l||l|l|}
		\cline{2-7}
		\multicolumn{1}{l|}{} & \multicolumn{2}{c||}{\textbf{Minimax}} & \multicolumn{2}{c||}{\textbf{Lexicographic Minimax}} & \multicolumn{2}{c|}{\textbf{Minimum Sum}} \\ \hline
		\multicolumn{1}{|c|}{Buses} & \multicolumn{1}{c|}{\begin{tabular}[c]{@{}c@{}}Average\\ Change\end{tabular}} & \multicolumn{1}{c||}{\begin{tabular}[c]{@{}c@{}}Standard\\ Deviation\end{tabular}} & \multicolumn{1}{c|}{\begin{tabular}[c]{@{}c@{}}Average\\ Change\end{tabular}} & \multicolumn{1}{c||}{\begin{tabular}[c]{@{}c@{}}Standard\\ Deviation\end{tabular}} & \multicolumn{1}{c|}{\begin{tabular}[c]{@{}c@{}}Average\\ Change\end{tabular}} & \multicolumn{1}{c|}{\begin{tabular}[c]{@{}c@{}}Standard \\ Deviation\end{tabular}} \\ \hline
		\multicolumn{1}{|l|}{$z^*$} & 30.77 & 26.70 & 26.00 & 24.35 & 24.47 & 25.04 \\ \hline
		\multicolumn{1}{|l|}{$z^*+1$} & 27.57 & 25.30 & 22.70 & 22.69 & 19.67 & 23.61 \\ \hline
		\multicolumn{1}{|l|}{$z^*+2$} & 25.57 & 23.61 & 20.80 & 21.29 & 18.20 & 22.46 \\ \hline
		\multicolumn{1}{|l|}{$z^*+3$} & 23.87 & 23.36 & 18.43 & 19.87 & 16.67 & 21.04 \\ \hline
		\multicolumn{1}{|l|}{$z^*+4$} & 22.37 & 22.09 & 19.10 & 18.91 & 15.23 & 19.67 \\ \hline
		\multicolumn{1}{|l|}{$z^*+5$} & 22.27 & 21.09 & 15.67 & 17.72 & 14.17 & 19.10 \\ \hline
		\multicolumn{1}{|l|}{$z^*+6$} & 19.23 & 20.78 & 14.90 & 17.05 & 13.20 & 17.94 \\ \hline
		\multicolumn{1}{|l|}{$z^*+7$} & 19.37 & 19.43 & 14.97 & 16.47 & 12.50 & 17.23 \\ \hline
		\multicolumn{1}{|l|}{$z^*+8$} & 17.23 & 18.52 & 13.77 & 15.38 & 11.43 & 16.54 \\ \hline
		\multicolumn{1}{|l|}{$z^*+9$} & 16.97 & 17.57 & 12.57 & 14.57 & 10.67 & 15.77 \\ \hline
		\multicolumn{1}{|l|}{$z^*+10$} & 15.30 & 16.64 & 11.10 & 14.28 & 9.90 & 15.03 \\ \hline
	\end{tabular}
	\caption{Average and standard deviation of absolute start time changes within each instance (in minutes), Scenario II}
\end{table}

\begin{table}[!h]
	\centering
	\begin{tabular}{c|l|l||l|l||l|l|}
		\cline{2-7}
		\multicolumn{1}{l|}{} & \multicolumn{2}{c||}{\textbf{Minimax}} & \multicolumn{2}{c||}{\textbf{Lexicographic Minimax}} & \multicolumn{2}{c|}{\textbf{Minimum Sum}} \\ \hline
		\multicolumn{1}{|c|}{Buses} & \multicolumn{1}{c|}{\begin{tabular}[c]{@{}c@{}}Average\\ Change\end{tabular}} & \multicolumn{1}{c||}{\begin{tabular}[c]{@{}c@{}}Standard\\ Deviation\end{tabular}} & \multicolumn{1}{c|}{\begin{tabular}[c]{@{}c@{}}Average\\ Change\end{tabular}} & \multicolumn{1}{c||}{\begin{tabular}[c]{@{}c@{}}Standard\\ Deviation\end{tabular}} & \multicolumn{1}{c|}{\begin{tabular}[c]{@{}c@{}}Average\\ Change\end{tabular}} & \multicolumn{1}{c|}{\begin{tabular}[c]{@{}c@{}}Standard \\ Deviation\end{tabular}} \\ \hline
		\multicolumn{1}{|l|}{$z^*$} & 22.57 & 18.31 & 19.17 & 17.28 & 16.90 & 18.74 \\ \hline
		\multicolumn{1}{|l|}{$z^*+1$} & 11.13 & 11.82 & 8.57 & 10.66 & 7.53 & 12.30 \\ \hline
		\multicolumn{1}{|l|}{$z^*+2$} & 6.90 & 8.78 & 5.90 & 8.01 & 4.97 & 8.52 \\ \hline
		\multicolumn{1}{|l|}{$z^*+3$} & 6.47 & 7.86 & 4.03 & 6.48 & 3.30 & 6.83 \\ \hline
		\multicolumn{1}{|l|}{$z^*+4$} & 4.93 & 6.36 & 3.13 & 5.40 & 2.50 & 6.08 \\ \hline
		\multicolumn{1}{|l|}{$z^*+5$} & 3.63 & 5.27 & 2.23 & 4.56 & 1.77 & 4.94 \\ \hline
		\multicolumn{1}{|l|}{$z^*+6$} & 3.27 & 4.12 & 1.63 & 3.63 & 1.47 & 4.21 \\ \hline
		\multicolumn{1}{|l|}{$z^*+7$} & 2.70 & 3.50 & 1.20 & 2.81 & 1.13 & 3.09 \\ \hline
		\multicolumn{1}{|l|}{$z^*+8$} & 1.53 & 2.31 & 0.70 & 1.88 & 0.70 & 2.07 \\ \hline
		\multicolumn{1}{|l|}{$z^*+9$} & 1.20 & 1.63 & 0.47 & 1.25 & 0.47 & 1.25 \\ \hline
		\multicolumn{1}{|l|}{$z^*+10$} & 0.77 & 1.24 & 0.27 & 1.00 & 0.27 & 1.00 \\ \hline
	\end{tabular}
	\caption{Average and standard deviation of absolute start time changes within each instance (in minutes), Scenario III}
\end{table}

\pagebreak
\begin{figure}[!hp]
	\includegraphics[width=\textwidth]{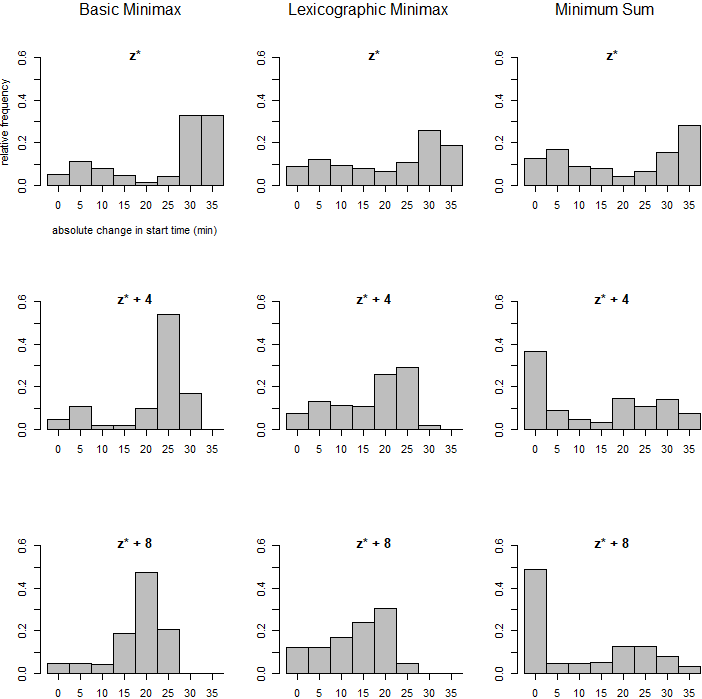}
	\centering 
	\caption{Start time change distribution, Scenario I}
\end{figure}

\begin{figure}[!hp]
	\centering
	\includegraphics[width=\textwidth]{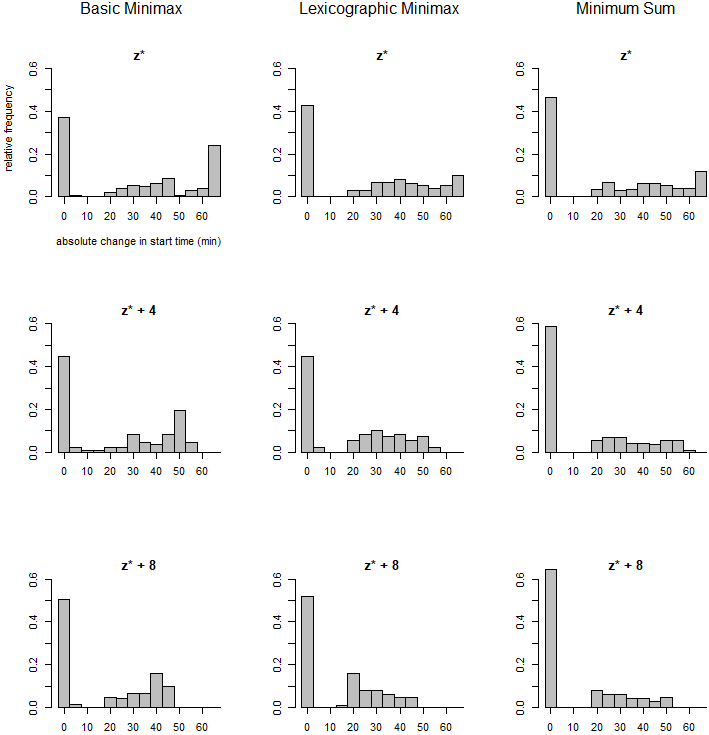} 
	\caption{Start time change distribution, Scenario II}
\end{figure} 

\begin{figure}[!ht]
	\centering
	\includegraphics[width=\textwidth]{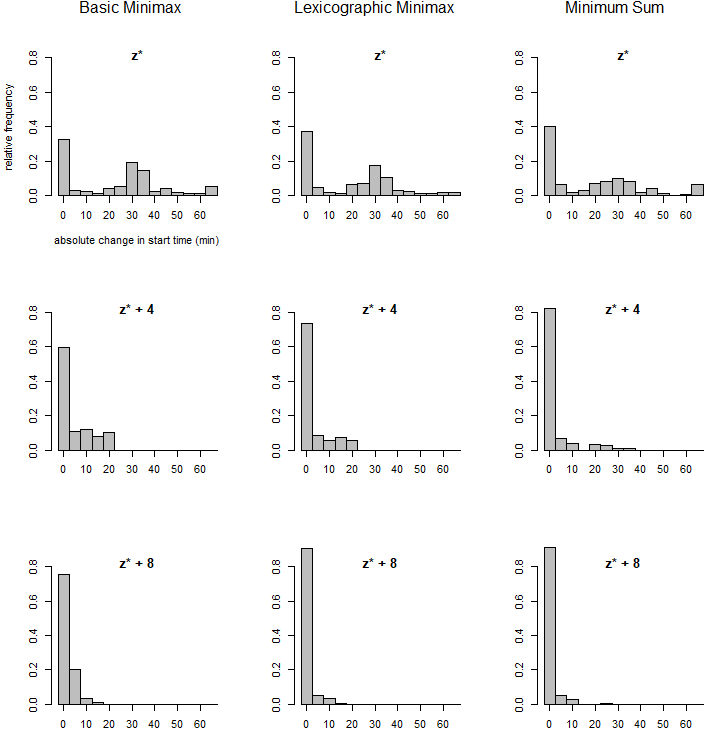}
	\caption{Start time change distribution, Scenario III}
\end{figure} 

\end{document}